\documentclass[a4paper, 12pt]{article}
\usepackage{CJKutf8}
\usepackage{amsmath,amssymb,amsthm,amscd,mathrsfs,amsfonts}
\usepackage{graphicx}
\usepackage{indentfirst}
\usepackage{url}
\usepackage{enumitem}
\usepackage[numbers,sort&compress]{natbib}
\date{}
\numberwithin{equation}{section}
\begin{document}
\begin{CJK}{UTF8}{gkai}
\newtheorem{thm}{Theorem}[section]
\newtheorem{pro}[thm]{Proposition}
\newtheorem{lem}[thm]{Lemma}
\newtheorem{remark}[thm]{Remark}
\newtheorem{cor}[thm]{Corollary}
\newtheorem{ques}[thm]{Question}
\theoremstyle{definition}
\newtheorem{definition}[thm]{Definition}
\newtheorem{ex}[thm]{Example}
\renewcommand{\thefootnote}{}
\baselineskip 16pt \setcounter{section}{0}

\title{On the weak norm of $\mathscr{U}_p$-residuals of all subgroups  of a finite group\thanks{Supported by the  major project of Basic and
Applied Research (Natural Science)  in Guangdong Province, China (Grant Number:  2017KZDXM058) and
the Science and Technology Program of Guangzhou Municipality, China (Grant number: 201804010088).}}
\author{Yubo Lv\thanks{School of Mathematical Sciences, Guizhou Normal University, Guiyang, 550001, China, \texttt{lvyubo341281@163.com}
}\and Yangming Li\thanks{Corresponding author}\ \thanks{College of Mathematics, Guizhou Normal University, Guiyang, 550001, China;
Dept.\ of Math., Guangdong University of Education,Guangzhou, 510310, China, \texttt{liyangming@gdei.edu.cn}}
}

\maketitle
\medskip
\noindent 
\begin{abstract}
  Let $\mathscr{F}$ be a formation and $G$ a finite group.  The weak norm of a subgroup $H$ in $G$ with respect to  $\mathscr{F}$ is defined by $N_{\mathscr{F}}(G,H)=\underset{T\leq H}{\bigcap}N_G(T^{\mathscr{F}})$. In particular,  $N_{\mathscr{F}}(G)=N_{\mathscr{F}}(G,G)$.  Let $N^i_{\mathscr{F}}(G)$,$i\geq 1$, be a upper series of $G$ by setting $N^0_{\mathscr{F}}(G)=1$, $N^{i+1}_{\mathscr{F}}(G)/N^i_{\mathscr{F}}(G)=N_{\mathscr{F}}(G/N^i_{\mathscr{F}}(G))$ and denoted by $N^{\infty}_{\mathscr{F}}(G)$ the terminal term of the series. In this paper, for  the case $\mathscr{F}\in\{\mathscr{U}_p,\mathscr{U}\}$, where $\mathscr{U}_p$($\mathscr{U}$,respectively) is the class of all finite $p$-supersolvable groups(supersolvable groups,respectively), we characterize the structure of some given finite groups  by the properties of  weak norm of some subgroups in $G$ with respect to $\mathscr{F}$.  Some of our main results may regard as a continuation of many nice previous work.
\end{abstract}


\section{Introduction}

Throughout this paper, all groups are finite. We use the standard terminology and notations as in \cite{Hup67}. For the specific, we denote

\begin{itemize}
  \item $G$: a finite group.
  \item $|G|$: the order of $G$.
  \item $\pi(G)$: the set of prime divisors of $|G|$.
  \item $p$: a prime.
  \item $G_p$: the Sylow $p$-subgroup of $G$.
  \item $F_p(G)$: the $p$-Fitting subgroup of $G$.
  \item $\mathfrak{G}$: the class of all finite groups.
  \item $\mathscr{F}$: a formation of groups, that is a class of finite groups satisfying the following:

  (1) if $G\in\mathscr{F}$ and $N$ is a normal subgroup of $G$, then $G/N\in\mathscr{F}$, and

  (2) if $N_1$,$N_2$ are normal subgroups of $G$ such that $G/N_i\in\mathscr{F}(i=1,2)$, then $G/(N_1\cap N_2)\in\mathscr{F}$.

  \item $G^{\mathscr{F}}$: the $\mathscr{F}$-residual of $G$, that is the intersection of all those normal subgroups $N$ of $G$ such that $G/N\in\mathscr{F}$.
  \item $\mathscr{F}_1\mathscr{F}_2$: the formation product or Gasch\"{e}tz product of $\mathscr{F}_1$ and $\mathscr{F}_2$, that is the class $\{G\in\mathfrak{G} \mid G^{\mathscr{F}_2}\in\mathscr{F}_1 \}$. In particular, denote $\mathscr{F}^2=\mathscr{F}\mathscr{F}$.
  \item $\mathscr{A}$: the class of all Abelian groups.
  \item $\mathscr{N}$: the class of all nilpotent groups.
  \item $\mathscr{N}_p$: the class of all $p$-nilpotent groups.
  \item $\mathscr{U}$: the class of all supersolvable groups.
  \item $\mathscr{U}_p$: the class of all $p$-supersolvable groups.
  \item $l_p(G)$: the $p$-length of a $p$-solvable group $G$, that is the number of the $p$-factor groups in the upper $p$-series of $G$:
$$1=P_0(G)\unlhd M_0(G)\unlhd P_1(G)\unlhd M_1(G)\unlhd \cdots \unlhd P_n(G)\unlhd M_n(G)=G$$ such that $M_i(G)/P_i(G)=O_{p'}(G/P_i(G))$ and $P_i(G)/M_{i-1}(G)=O_p(G/M_{i-1}(G))$.

  \item $h_p(G)$: the $p$-Fitting length of a $p$-solvable group $G$, that is the  smallest positive integer $n$ such that $$1=F^0_p(G)\leq F^1_p(G)\leq \cdots\leq F^{n-1}_p(G)\leq F^n_p(G)=G,$$ where $F^1_p(G)=F_p(G)$ and $F^{i+1}_p(G)/F^i_p(G)=F_p(G/F^i_p(G))$ for $i=1,2,\cdots,n-1$.
  \item $h(G)$: the Fitting  length of a solvable group $G$, that is  the positive integer $n$ such that $1=F^0(G)\leq F^1(G)\leq \cdots\leq F^{n-1}(G)\leq F^n(G)=G,$ where $F^1(G)=F(G)$ and $F^{i+1}(G)/F^i(G)=F(G/F^i(G))$ for $i=1,2,\cdots,n-1$.
  \item $N(G)$: the intersection of the normalizers of all subgroups of $G$.
  \item $S(G)$(or $N_{\mathscr{N}}(G)$):  the intersection of the normalizers  of the nilpotent residuals of all subgroups of $G$.
  \item $N^{\mathscr{N}_p}(G)$: the intersection of the normalizers  of the $p$-nilpotent residuals of all subgroups of $G$.

\end{itemize}

It is interesting to characterize the structure of a given group by using of some special subgroups. For example, Gash\"{e}tz and N. It\^{o}                  \cite[Satz 5.7,p.436]{Hup67} proved  that $G$ is solvable with Fitting length at most 3 if all minimal subgroups of $G$ are normal.  Also, it is well know that $G$ is nilpotent if $G'$ normalizers each subgroup of $G$(see Baer's theorem in \cite{Baer56}). Further more,  R. Baer in \cite{Baer35} defined the subgroup $N(G)$, the norm of a group $G$. Obviously, a group $G$ is a Dedeking group if and only if $G=N(G)$. The norm of a group has many other good properties and has studied further by many scholars. In recent years, some weaker versions of the concepts of norm of groups  have been introduced. Let $\mathscr{F}$ be a non-empty formation. Recently, Su and Wang in \cite{SuNing13,SuNing14} introduced the subgroup $N_{\mathscr{F}}(G)$, the norm of $\mathscr{F}$-residual of a group $G$ as follows

\begin{equation*}
N_{\mathscr{F}}(G)= \underset{H\leq G}{\bigcap}{N_G(H^{\mathscr{F}})}. \tag{$\star$}
\end{equation*}

As in \cite{SuNing13}, we set $N^0_{\mathscr{F}}(G)=1$ and if $N^i_{\mathscr{F}}(G)$ is defined, set $N^{i+1}_{\mathscr{F}}(G)/N^i_{\mathscr{F}}(G)=N_{\mathscr{F}}(G/N^i_{\mathscr{F}}(G))$. The subgroup $N^{\infty}_{\mathscr{F}}(G)$ is the terminal term of the ascending series. In fact, $N^{\infty}_{\mathscr{F}}(G)=N^k_{\mathscr{F}}(G)$ for some integer $k$ such that $N^k_{\mathscr{F}}(G)=N^{k+1}_{\mathscr{F}}(G)=N^{k+2}_{\mathscr{F}}(G)=\cdots$.

Many scholars also call $N_{\mathscr{F}}(G)$($\mathscr{F}=\mathscr{A},\mathscr{N}, \mathscr{N}_p$, respectively)  the {\em generalized norm} of   group $G$. Obviously, $N_{\mathscr{F}}(G)$ is a characteristic subgroup and every element of $N_{\mathscr{F}}(G)$ normalize the $\mathscr{F}$-residual of each subgroup of $G$.  The so called norm of $\mathscr{F}$-residual has many other nice properties and also closely related to the  global properties of a given group. For some given formation $\mathscr{F}$, there are many papers devoted to study the $p$-length, Fitting length, solvability, ($p$-)nilpotency  and so on. For example,  for the case $\mathscr{F}=\mathscr{A}$, Li and  Shen in \cite{LiSR10} denoted $N_{\mathscr{\mathscr{A}}}(G)$ by $D(G)$. They  fingered out that $G$ is solvable with Fitting length at most 3 if all elements of $G$ of prime order are in $D(G)$(see \cite[Theorem 4.1]{LiSR10}). It is a dual problem of Gash\"{e}tz and N. It\^{o}\cite[Satz 5.7,p.436]{Hup67}. Li and Shen also defined the $D$-group, i.e., $G=D(G)$, they characterized the relationship between $D(G)$ and $G$. Shen, Shi and Qian in \cite{ShenZC12} considered the case  $\mathscr{F}=\mathscr{N}$, they denoted $N_{\mathscr{N}}(G)$ by $S(G)$. Shen et al., deeply studied the dual problem of Gash\"{e}tz and N. It\^{o} and characterized the  $\mathscr{F}_{nn}$-groups by means of the subgroup $S^{\infty}(G)$, where $S^{\infty}(G)=N^{\infty}_{\mathscr{N}}(G)$ and $\mathscr{F}_{nn}$-groups are class of groups belong to $\mathscr{N}\mathscr{N}$. They also introduced the $S$-group, i.e., $G=S(G)$ and given some sufficient and necessary conditions involved $S$-groups. Meanwhile, Gong and Guo in \cite{GongL13} also consider the case  $\mathscr{F}=\mathscr{N}$ and given some meaningful conclusions. In the case $\mathscr{F}=\mathscr{N}_p$, Guo and Li in \cite{LiXuan15} introduced the norm of $\mathscr{N}_p$-residual of a group $G$, they denoted $N_{\mathscr{N}_p}(G)$ by $N^{\mathscr{N}_p}(G)$. As a local version of Gong's results, Li and Guo characterized the relationship between $C_G(G^{\mathscr{N}_p})$ and $N^{\mathscr{N}_p}(G)$. In particular, Li and Guo in \cite{LiXuan15} also investigated the relationship between $N_{\mathscr{N}}(G)$ and $N^{\mathscr{N}_p}(G)$.  For more detail and other relevant conclusions about the norm of $\mathscr{F}$-residual, please see  \cite{ShenZC12,SuNing13,LiSR10,SuNing14,LiXuan16,LiXuan15,ShenZC14,WangJX07,GongL13,HuB20}.

We wonder whether the above conclusions hold for general formations.  A natural idea is to replace $\mathscr{F}(\mathscr{F}\in\{\mathscr{N}_p,\mathscr{N},\mathscr{A}\})$ by $\mathscr{U}_p$ or $\mathscr{U}$ in $(\star)$.

\begin{remark}

(1) In general, for a group $G$ and some $p\in\pi(G)$, $N_{\mathscr{U}_p}(G)\neq N_{\mathscr{F}}(G)$ is possible, where $\mathscr{F}\in\{\mathscr{N}_p,\mathscr{N},\mathscr{A}\}$. For example, let $G=S_4$, the symmetric group of degree $4$. Obviously, $G$ is a $3$-supersolvable non-$3$-nilpotent group, so $G=N_{\mathscr{U}_3}(G)$. Pick a subgroup $H\cong S_3$, the symmetric group of degree $3$, but $N_G(H^{\mathscr{N}_3})=N_G(S^{\mathscr{N}_3}_3)=N_G(C_3)=S_3<S_4$, so $N_{\mathscr{N}_3}(G)<G$.

(2) The condition in  definition $(\star)$  that the intersection of the $\mathscr{F}$-residuals  of all subgroups of $G$ may be too strong and some of the subgroups may be redundant whenever $\mathscr{F}=\mathscr{U}_p$(see Example \ref{exre} below).

(3) The case that $N_{\mathscr{U}}(G)=N_{\mathscr{U}_p}(G)=1$ is possible for a solvable group $G$ and $p\in\pi(G)$. Let $H=\langle a,b\mid a^3=b^3=1,[a,b]=1\rangle\cong C_3\times C_3$ and $Q_8=\langle c,d\mid c^4=1,c^2=d^2=e,c^d=c^{-1}\rangle$. Considering the irreducible action of $Q_8$  on $H$ by $a^c=a^{-1}b,b^c=ab,a^d=b^{-1},b^d=a$, denote $T=H\rtimes Q_8\cong (C_3\times C_3)\rtimes Q_8$. Let $C=\langle f\rangle$ be a cyclic group of order $3$ and let $C$ act on $T$ by $a^f=b^{-1},b^f=ab^{-1},c^f=d^3,d^f=cd$. Then $G=\langle a,b,c,d,e,f\rangle\cong ((C_3\times C_3)\rtimes Q_8)\rtimes C_3$ is solvable(IdGroup=[216,153]). It is easy to see that $G^{\mathscr{U}_2}=\langle a,b,c,d,e\rangle\cong (C_3\times C_3)\rtimes Q_8$ and $C_G(G^{\mathscr{U}_2})=1$.  By  Theorem \ref{thmCGNUp1}, $N_{\mathscr{U}_2}(G)=1$. Further more,  by Theorem \ref{lemNEU},  $N_{\mathscr{U}}(G)=1$. So  we have  $N_{\mathscr{U}}(G)=N_{\mathscr{U}_2}(G)=1$.

\end{remark}

Recall that the weak centralizer of $H$ in $G$, $C^{\ast}_G(H)$, introduced in \cite[P.33,Definitions]{mW}, is defined by $$C^{\ast}_G(H)=\bigcap\{N_G(K):K\leq H\}.$$ In the above investigation, we introduce the following  more general and interesting definition.

\begin{definition}
Let $G$ be a group and $\mathscr{F}$ a formation. We define $N_{\mathscr{F}}(G,H)$, the weak norm  of $H$ in $G$ with respect to  $\mathscr{F}$ as follows: $$N_{\mathscr{F}}(G,H)=\underset{T\leq H}{\bigcap}{N_G(T^{\mathscr{F}})}.$$

In particular, $N_{\mathscr{U}}(G,H)=\underset{T\leq H}{\bigcap}{N_G(T^{\mathscr{U}})}$ and $N_{\mathscr{U}_p}(G,H)=\underset{T\leq H}{\bigcap}{N_G(T^{\mathscr{U}_p})}$.
\end{definition}

\begin{ex}\label{exre}
Let $G=\langle a,b,c,d\rangle\rtimes\langle e,f\rangle\cong C^4_2\rtimes C_6$, where $e^2=f^3=1$ and $a^e=ac,b^e=bd,c^e=c,d^e=d,a^f=b,b^f=ab,c^f=d,d^f=cd$ (IdGroup=[96,70] in GAP \cite{GAP}). Let $H=\langle a,b,f\rangle$, then $N_{\mathscr{U}_2}(G)=N_{\mathscr{U}_2}(G,H)=\langle a,b,c,d,f\rangle\cong C^4_2\rtimes C_3<G$.
\end{ex}

Obviously, $N_{\mathscr{F}}(G)=N_{\mathscr{F}}(G,G)\leq N_{\mathscr{F}}(G,H)\leq G$. Without causing confusion, we call $N_{\mathscr{F}}(G)=N_{\mathscr{F}}(G,G)$ the norm of $G$ with respect to $\mathscr{F}$.  Moreover, for a subgroup $H$ of $G$,  we call $G$ is an $\mathscr{N}_1$-group  with respect to $H$ and $\mathscr{F}$ if $G=N_{\mathscr{F}}(G,H)$. In this paper, we mainly investigate  properties of $N_{\mathscr{F}}(G,H)$ and the influence of $N_{\mathscr{F}}(G,H)$ on the structure of group $G$.  Actually, we mainly consider the case $\mathscr{F}\in\{\mathscr{U},\mathscr{U}_p\}$. Our main work may be regard as the continuation of some conclusions in \cite{LiSR10,ShenZC12,LiXuan15,GongL13}.

\section{Preliminaries}
In this section, we always assume that $\mathscr{F}$ is a non-empty formation and $G$ is a group. We first give some important lemmas.

\begin{lem}\label{lemwn}
Let $H,K,N$ be subgroups of $G$ and $N\unlhd G$. Then

(1) If $H\leq K$, then $N_{\mathscr{F}}(G,K)\leq N_{\mathscr{F}}(G,H)$;

(2) $K\cap N_{\mathscr{F}}(G,H)\leq N_{\mathscr{F}}(K,K\cap H)$,  in particular, if $H\leq K$, then $K\cap N_{\mathscr{F}}(G,H)\leq N_{\mathscr{F}}(K,H)$.

(3) If $N\leq H$, then $N_{\mathscr{F}}(G,H)N/N\leq N_{\mathscr{F}}(G/N,H/N)$.
\end{lem}
\proof (1) By definition, $N_{\mathscr{F}}(G,K)=\underset{T\leq K}{\bigcap}{N_G(T^{\mathscr{F}})}\leq \underset{T\leq H}{\bigcap}{N_G(T^{\mathscr{F}})}=N_{\mathscr{F}}(G,H)$.

(2) Obviously, $K\cap N_{\mathscr{F}}(G,H)=K\cap (\underset{T\leq H}{\bigcap}{N_G(T^{\mathscr{F}})})\leq \underset{T\leq H\cap K}{\bigcap}{N_K(T^{\mathscr{F}})}=N_{\mathscr{F}}(K,H\cap K)$.  In particular, if $H\leq K$, then $K\cap N_{\mathscr{F}}(G,H)\leq N_{\mathscr{F}}(K,H)$.

(3) Let $x\in N_{\mathscr{F}}(G,H)$, then $x$ normalize $T^{\mathscr{F}}$ for every $T\leq H$, so $xN$ normalize $T^{\mathscr{F}}N/N=(TN/N)^{\mathscr{F}}$. Thus every element of $N_{\mathscr{F}}(G,H)N/N$ normalize $(T/N)^{\mathscr{F}}$ for all subgroups $T/N$ of $H/N$, so $N_{\mathscr{F}}(G,H)N/N\leq N_{\mathscr{F}}(G/N,H/N)$. \qed\\

As a corollary, we have

\begin{lem}\label{lemSQ}
Let $G$ be a group, let $K$ be a subgroup of $G$ and $N$ a  normal subgroup of $G$, then

(1) $K\cap N_{\mathscr{F}}(G)\leq N_{\mathscr{F}}(K)$.

(2) $N_{\mathscr{F}}(G)N/N\leq N_{\mathscr{F}}(G/N)$.
\end{lem}

\begin{lem}\cite[Lemmas 2.2,2.3]{SuNing13}\label{lemDnUp}
Let $G$ be a group, let $K$ be a subgroup of $G$ and $N$ a  normal subgroup of $G$, then

(1) $K\cap N^{\infty}_{\mathscr{F}}(G)\leq N^{\infty}_{\mathscr{F}}(K)$.

(2) $N^{\infty}_{\mathscr{F}}(G)N/N\leq N^{\infty}_{\mathscr{F}}(G/N)$.

(3) If $N\leq N^{\infty}_{\mathscr{F}}(G)$, then $N^{\infty}_{\mathscr{F}}(G)/N=N^{\infty}_{\mathscr{F}}(G/N)$.
\end{lem}

\begin{lem}\label{lempL}
Let $G$ be a $p$-solvable group.

(1) If $N\unlhd G$, then $l_p(G/N)\leq l_p(G)$.

(2) If $U\leq G$, then $l_p(U)\leq l_p(G)$.

(3) Let $N_1$ and $N_2$ be two normal subgroups of $G$, then $l_p(G/(N_1\cap N_2))\leq max\{l_p(G/N_1),l_p(G/N_2)\}$.

(4) $l_p(G/\Phi(G))=l_p(G)$.

(5) If $N$ is a normal $p'$-group of $G$, then $l_p(G/N)=l_p(G)$.
\end{lem}
\proof The proof of (1)-(4) follows from \cite[VI,6.4]{Hup67} and (5) is obviously. \qed

\begin{lem}\label{lempnPhi}
Let $G$ be a group and $N$ a normal subgroup of $G$. Then

(1) if $N/N\cap\Phi(G)$  is $p$-nilpotent, then $N$ is $p$-nilpotent;

(2) Let $H$ be a subgroup of $G$ and $N$ be a $p'$-group, if $HN/N$ is $p$-nilpotent, then $H$ is $p$-nilpotent.
\end{lem}
\proof The statement (1) is directly form  \cite[Lemma 2.5]{Skiba15}, and (2) is easy. \qed

\begin{lem}\label{lemhp}
Let $G$ be a $p$-solvable group and $H$ a subgroup of $G$, let $H,N,A,B$ be   subgroups of $G$ and $N,A,B$ are normal in $G$, then

(1) $h_p(H)\leq h_p(G)$.

(2) $h_p(G/N)\leq h_p(G)$.

(3) If $G=A\times B$, then $h_p(G)=\max\{h_p(A),h_p(B)\}$.

(4) If $h_p(G/A)\leq k$ and $h_p(G/N)\leq k$, then $h_p(G/(A\cap B))\leq k$.

(5) $h_p(G/\Phi(G)))=h_p(G)$.

(6) If $N$ is a $p'$-group, then $h_p(G/N)=h_p(G)$.
\end{lem}
\proof  Let $1=N_0\leq N_1\leq N_2\leq \cdots\leq N_r=G$ be the shortest normal chain of $G$ with $p$-nilpotent factors $N_i/N_{i-1}$ for all $i=1,2,\cdots,r$.

(1) If $H\leq G$, then $1=N_0\cap H\leq N_1\cap H\leq N_2\cap H\leq \cdots\leq N_r\cap H=G\cap H=H$ is a normal chain of $H$ with $N_i\cap H/N_{i-1}\cap H\cong (N_i\cap H)N_{i-1}/N_{i-1}$ a $p$-nilpotent factors for all $i=1,2,\cdots,r$, so $h_p(H)\leq h_p(G)$.

(2) If $N\unlhd G$, obviously, $\bar{1}=N_0N/N\leq N_1N/N\leq N_2N/N\leq \cdots\leq N_r/N=G/N$ is a  normal chain of $G/N$. Note that  $N_iN/N/N_{i-1}N/N=N_iN/N_{i-1}N\cong N_i/N_{i-1}(N_i\cap N)\leq N_i/N_{i-1}$ is $p$-nilpotent, so $h_p(G/N)\leq h_p(G)$.

(3) Let $1=A_0\leq A_1\leq A_2\leq \cdots\leq A_r=A$ and $1=B_0\leq B_1\leq B_2\leq \cdots\leq B_t=B$ be the shortest normal chain of $A$ and $B$  with $p$-nilpotent factors respectively. Without loss of generality, assume $r\leq t$, then $1=A_0B_0\leq A_1B_1\leq A_2B_2\leq \cdots\leq A_rB_r\leq A_rB_{r+1}\leq A_rB_{r+2}\leq A_rB_t=AB$ is a normal chain. Since $A_iB_i/A_{i-1}B_{i-1}\cong A_i/A_{i-1}(A_i\cap B_{i-1})\cdot B_i/B_{i-1}(B_i\cap A_{i-1})$  whenever $i\leq r$ and $A_rB_j/A_rB_{j-1}\cong B_j/B_{j-1}(B_j\cap A_i)$  whenever $r<j\leq t$ are   $p$-nilpotent. So $h_p(G)=\max\{h_p(A),h_p(B)\}$.

(4) Since $G/(A\cap B)\cong G/A\times G/B$, the result  follows from (1) and (3).

(5) Assume that $\bar{1}=\Phi(G)/\Phi(G)=T_0/\Phi(G)\leq T_1/\Phi(G)\leq T_2/\Phi(G)\leq \cdots\leq T_s/\Phi(G)=G/\Phi(G)$ is the normal chain of $G$ with $p$-nilpotent factors, so $s\leq r$ and $T_i/T_{i-1}\cong T_i/\Phi(G)/T_{i-1}/\Phi(G)$ is $p$-nilpotent.  Since $T_1/\Phi(G)$ is $p$-nilpotent, then  $T_1$ is $p$-nilpotent by Lemma \ref{lempnPhi}, so $1\leq T_1\leq T_2\leq \cdots \leq T_s=G$ is a normal chain with $p$-nilpotent factors. Now we have $r\leq s$, hence $s=r$.

(6) Let $\bar{1}=N_0/N\leq N_1/N\leq N_2/N\leq \cdots\leq N_r/N=G/N$ be a normal chain of $G/N$ with $p$-nilpotent factors $N_i/N/N_{i-1}/N\cong N_i/N_{i-1}$ for all $i=1,2,\cdots,r$. Now by Lemma \ref{lempnPhi}, $N_1$ is $p$-nilpotent  since $N$ is a $p'$-group, so $1\leq N_1\leq N_2\leq \cdots\leq N_r=G$ is a normal chain of $G$ with $p$-nilpotent factors $N_i/N_{i-1}$ for all $i=1,2,\cdots,r$, which implies that $h_p(G)\leq h_p(G/N)$, now by (1), we have $h_p(G/N)=h_p(G)$. \qed\\

As a local vision of \cite[Lemma 3.2]{Jabara18}, we have
\begin{lem}\label{lempl21}
Let $G$ be a solvable group, then $l_p(G)\leq 1$ if $h(G)\leq 2$ for every $p\in\pi(G)$.
\end{lem}

\begin{lem}\cite[Theorem 1 and Proposition 1]{Adolf96}\label{lemAdolf}
Let $\mathscr{F}$ be a saturated formation.

(1) Assume that $G$ is a group such that $G$ does not belong to $\mathscr{F}$, but all its proper subgroups belong to $\mathscr{F}$. Then $F'(G)/\Phi(G)$  is the unique minimal normal subgroup of $G/\Phi(G)$, where $F'(G)=Soc(G \mod \Phi(G))$, and $F'(G)=G^{\mathscr{F}}\Phi(G)$. In addition, if the derived subgroup of $G^{\mathscr{F}}$ is a proper subgroup of $G^{\mathscr{F}}$, then $G^{\mathscr{F}}$ is a soluble group. Furthermore, if $G^{\mathscr{F}}$ is soluble, then  $F'(G)=F(G)$, the Fitting subgroup of $G$. Moreover $(G^{\mathscr{F}})'=T\cap G^{\mathscr{F}}$ for every maximal subgroup $T$ of $G$ such that $G/T_G\not\in \mathscr{F}$  and $F'(G)T=G$.

(2) Assume that $G$ is a group such that $G$ does not belong to $\mathscr{F}$ and there exists a maximal subgroup $M$ of $G$ such that $M\in\mathscr{F}$ and  $G=MF(G)$. Then $G^{\mathscr{F}}/(G^{\mathscr{F}})'$ is a chief factor of $G$, $G^{\mathscr{F}}$ is a $p$-group for some prime $p$, $G^{\mathscr{F}}$ has exponent $p$ if $p>2$ and exponent at most $4$ if $p=2$. Moreover, either $G^{\mathscr{F}}$ is elementary abelian or $(G^{\mathscr{F}})' = Z(G^{\mathscr{F}})=\Phi(G^{\mathscr{F}})$ is an elementary abelian group.
\end{lem}

\begin{lem}\label{lemGUpp}\label{lemHUpp}
Let $G$ be a  $p$-solvable minimal  non-$p$-supersolvable group, then $G^{\mathscr{U}_p}$ is a $p$-group.
\end{lem}
\proof  By Lemma \ref{lemAdolf}(1), $G/\Phi(G)$ has the unique minimal normal subgroup $F'(G)/\Phi(G)$, where $F'(G)=G^{\mathscr{U}_p}\Phi(G)$, so $F'(G)/\Phi(G)$ is an elementary abelian $p$-group since $G$ is a $p$-solvable group. Otherwise, $F'(G)/\Phi(G)$ is a $p'$-group and the $p$-supersolvability of $G/\Phi(G)/F'(G)/\Phi(G)$ implies that $G/\Phi(G)$ is $p$-supersolvable, thus $G$ is $p$-supersolvable, a contradiction. Furthermore, $\Phi(G)=1$. Now we have $G^{\mathscr{U}_p}$ is  solvable and $F'(G)=F(G)$, $G=F(G)T$ for some maximal subgroup $T$ of $G$. Now by Lemma \ref{lemAdolf}(2), $G^{\mathscr{U}_p}$ is a $p$-group. \qed

\begin{lem}\label{lemGUUp}
Let $G$ be a group and $\pi(G)=\{p_1,p_2,\cdots,p_n\}$, then $G^{\mathscr{U}}=G^{\mathscr{U}_{p_1}}G^{\mathscr{U}_{p_2}}\cdots G^{\mathscr{U}_{p_n}}$.
\end{lem}
\proof It follows from the supersolvability of $G/G^{\mathscr{U}}$ that $G^{\mathscr{U}_{p_i}}\leq G^{\mathscr{U}}$ for all $p_i\in\pi(G)$, so $G^{\mathscr{U}_{p_1}}G^{\mathscr{U}_{p_2}}\cdots G^{\mathscr{U}_{p_n}}\leq G^{\mathscr{U}}$. Conversely, obviously, $G/(G^{\mathscr{U}_{p_1}}G^{\mathscr{U}_{p_2}}\cdots G^{\mathscr{U}_{p_n}})\in\mathscr{U}_{p_i}$ for any $p_i\in\pi(G)$, so $G^{\mathscr{U}}\leq G^{\mathscr{U}_{p_1}}G^{\mathscr{U}_{p_2}}\cdots G^{\mathscr{U}_{p_n}}$, as desired. \qed

\section{Basic properties}

Now we give some elementary properties. Firstly, we have

\begin{pro}\label{lemNGAB}
Let $G=A\times B$ with $(|A|,|B|)=1$ and $H\leq G$,  then $N_{\mathscr{U}_p}(G,H)=N_{\mathscr{U}_p}(A,A\cap H)\times B$ whenever $p\in\pi(A)$ or $A\times N_{\mathscr{U}_p}(B,B\cap H)$ whenever $p\in\pi(B)$.
\end{pro}
\proof  By hypothesis, we may assume that $p\in\pi(A)\setminus\pi(B)$. Note that $H=(H\cap A)\times (H\cap B)$ and $T=(T\cap A)\times (T\cap B)$, so $T^{\mathscr{U}_p}=(T\cap A)^{\mathscr{U}_p}\times (T\cap B)^{\mathscr{U}_p}=(T\cap A)^{\mathscr{U}_p},$ where $T\leq H$. Now $N_G(T^{\mathscr{U}_p})=N_{A\times B}(T^{\mathscr{U}_p})=N_A((T\cap A)^{\mathscr{U}_p})\times B.$ We have
\begin{align*}
  N_{\mathscr{U}_p}(G,H) &=\underset{T\leq H}{\bigcap}{N_G(T^{\mathscr{U}_p})}  \notag \\
   {}&=(\underset{T\leq H}{\bigcap}{N_A((T\cap A)^{\mathscr{U}_p})})\times B  \notag\\
   {}&=(\underset{T\leq H\cap A}{\bigcap}{N_A(T^{\mathscr{U}_p})})\times B \\
   {}&=N_{\mathscr{U}_p}(A,A\cap H)\times B.
\end{align*}

Similarly, if  $p\in\pi(B)\backslash\pi(A)$, we have $N_{\mathscr{U}_p}(G,H)=A\times N_{\mathscr{U}_p}(B,B\cap H)$. \qed

\begin{pro}\label{proN1}
Let $G$ be a $p$-solvable group, then the $\mathscr{U}_p$-residual of $N_{\mathscr{U}_p}(G)$ is a $p$-group.
\end{pro}
\proof Denote $X=N_{\mathscr{U}_p}(G)$. We first prove this result in the case that $X=G$.  Assume $G\in\mathscr{U}_p$, obviously $G^{\mathscr{U}_p}=1$ is a $p$-group.   Assume now that $G\not\in\mathscr{U}_p$, if $O_p(G)>1$, then by Lemma \ref{lemSQ}(2), $G/O_p(G)=N_{\mathscr{U}_p}(G/O_p(G))$. By induction on $|G|$,  $(G/O_p(G))^{\mathscr{U}_p}=N^{\mathscr{U}_p}_{\mathscr{U}_p}(G/O_p(G))$ is a $p$-group. Obviously $G^{\mathscr{U}_p}$ is a $p$-group. If $O_p(G)=1$, then for any $K<G$, by Lemma \ref{lemSQ}(1), $K=N_{\mathscr{U}_p}(K)$, so by induction, $K^{\mathscr{U}_p}$ is a $p$-group. Note that $K^{\mathscr{U}_p}\unlhd G$, so $K^{\mathscr{U}_p}\leq O_p(G)=1$, this implies that $G$ is a $p$-solvable minimal non-$p$-supersolvable group, now by Lemma \ref{lemHUpp}, $G^{\mathscr{U}_p}$ is a $p$-group.

Assume now that $X<G$, by Lemma \ref{lemSQ}(1), $X=X\cap N_{\mathscr{U}_p}(G)\leq N_{\mathscr{U}_p}(X)\leq X$, so $X=N_{\mathscr{U}_p}(X)$. As a similar argument above, we also have the conclusion. \qed\\

Denote $Z_{\infty}(G)$ be the terminal term of the ascending central series of $G$.    As we know, $Z_{\infty}(G)=\bigcap\{N\unlhd G \mid  Z(G/N)=1\}$, we have the following similar result.

\begin{thm}
Let $G$ be a group, then $N^{\infty}_{\mathscr{U}_p}(G)=\bigcap\{N\unlhd G\mid N_{\mathscr{U}_p}(G/N)=1\}$.
\end{thm}
\proof Let $N\unlhd G$ such that $N_{\mathscr{U}_p}(G/N)=1$, then $N^{\infty}_{\mathscr{U}_p}(G/N)=1$. By Lemma \ref{lemDnUp}(2), $N^{\infty}_{\mathscr{U}_p}(G)N/N\leq N^{\infty}_{\mathscr{U}_p}(G/N)=1$, so $N^{\infty}_{\mathscr{U}_p}(G)\leq N$, thus $N^{\infty}_{\mathscr{U}_p}(G)\leq \bigcap\{N\unlhd G\mid N_{\mathscr{U}_p}(G/N)=1\}$.

Conversely, choice the least positive integer $n$ such that $N^{\infty}_{\mathscr{U}_p}(G)=N^n_{\mathscr{U}_p}(G)$. By definition, $N_{\mathscr{U}_p}(G/N^n_{\mathscr{U}_p}(G))=N^{n+1}_{\mathscr{U}_p}(G)/N^n_{\mathscr{U}_p}(G)=N^{n}_{\mathscr{U}_p}(G)/N^n_{\mathscr{U}_p}(G)=1$. Obviously, $\bigcap\{N\unlhd G\mid N_{\mathscr{U}_p}(G/N)=1\}\leq N^{\infty}_{\mathscr{U}_p}(G)$. \qed\\

Now we give some characterizations  between $N_{\mathscr{U}_p}(G,H)$ and $C_G(H^{\mathscr{U}_p})$.

\begin{thm}\label{thmCGNUp1}
Let $G$ be a $p$-solvable group and $H$ a normal subgroup of $G$, then $N_{\mathscr{U}_p}(G,H)=1$ if and only if $C_G(H^{\mathscr{U}_p})=1$.
\end{thm}
\proof  As we know, $C_G(H^{\mathscr{U}_p})\leq N_{\mathscr{U}_p}(G,H)$, the necessary is obviously.  Assume now $C_G(H^{\mathscr{U}_p})=1$, we prove that $N_{\mathscr{U}_p}(G,H)=1$. If not, let $N$ be a minimal normal subgroup of $G$ contained in $N_{\mathscr{U}_p}(G,H)$, if $N\nleq H$, then $N$ normalize $H^{\mathscr{U}_p}$ by definition, we have $[N,H^{\mathscr{U}_p}]\leq N\cap H^{\mathscr{U}_p}=1$ by the minimal normality of $N$, so $N\leq C_G(H^{\mathscr{U}_p})=1$, a contradiction. If $H<G$, then $(H,H)$ satisfies our hypothesis. By induction  on $|G||H|$ and Lemma \ref{lemwn}(2), $N\leq H\cap N_{\mathscr{U}_p}(G,H)\leq N_{\mathscr{U}_p}(H,H)=1$, a contradiction, thus $H=G$.  Now by the $p$-solvability of $G$, $N$ is a $p'$-group or an elementary abelian $p$-group. If $N$ is a $p'$-group, then $G/C_G(N)\in\mathscr{U}_p$, otherwise, there is  a minimal non-$p$-supersolvable subgroup $T/C_G(N)$ of $G/C_G(N)$. Let $T=C_G(N)L$ such that $C_G(N)\cap L\leq\Phi(L)$. Since $L/\Phi(L)\leq T/C_G(N)\cong L/L\cap C_G(N)$ is a minimal non-$p$-supersolvable group, we have $(L/\Phi(L))^{\mathscr{U}_p}=L^{\mathscr{U}_p}\Phi(L)/\Phi(L)$ is a $p$-group by Lemma \ref{lemGUpp}. Now let $P\in Syl_p(L^{\mathscr{U}_p})$, then $L^{\mathscr{U}_p}=P(L^{\mathscr{U}_p}\cap\Phi(L))$. By Frattini argument, $L=N_L(P)L^{\mathscr{U}_p}=N_L(P)P(L^{\mathscr{U}_p}\cap\Phi(L))=N_L(P)$, we have $P\unlhd L^{\mathscr{U}_p}$ and $N$ normalize $P$, so $P\leq C_G(N)\cap L\leq \Phi(L)$, thus $(L/\Phi(L))^{\mathscr{U}_p}=1$, $L\in\mathscr{U}_p$, contrary to our choice of $T/C_G(N)$. Consequently, $G^{\mathscr{U}_p}\leq C_G(N)$, hence $N\leq C_G(G^{\mathscr{U}_p})=1$, a contradiction. So $N$ is an elementary abelian $p$-group and $N\leq C_G(N)$. Note that $\Phi(G)\leq F(G)\leq C_G(N)$, we consider the following two cases:

\textbf{Case 1}. $\Phi(G)=C_G(N)$.

In this case, $F(G)=C_G(N)=\Phi(G)\leq F_p(G)$, so $N\leq\Phi(G)$. By the $p$-solvability of $G$, $F_p(G/\Phi(G))=F_p(G)/\Phi(G)\neq 1$, so $\Phi(G)=F(G)<F_p(G)$. Now follows from   $F_p(G)/O_{p'}(G)=O_{p'p}(G)/O_{p'}(G)=O_p(G/O_{p'}(G))$ that $O_{p'}(G)>1$, so $O_{p'}(G)\leq C_G(N)=\Phi(G)$. On the other hand, Let $P\in Syl_p(F_p(G))$, then $F_p(G)=[O_{p'}(G)]P$. By the Frattini argument, we have $G=F_p(G)N_G(P)=O_{p'}(G)N_G(P)=N_G(P)$, so $F_p(G)=P\times O_{p'}(G)=F(G)$, a contradiction.

\textbf{Case 2}. $\Phi(G)<C_G(N)$.

In this case, if $N\nleq \Phi(G)$, then there exists some maximal subgroup $M$ of $G$ such that $G=NM$ and $N\cap M=1$. Obviously, $(G/N)^{\mathscr{U}_p}=G^{\mathscr{U}_p}N/N\cong G^{\mathscr{U}_p}/G^{\mathscr{U}_p}\cap N\cong  M^{\mathscr{U}_p}$ and $N$ normalize $M^{\mathscr{U}_p}$. By the minimal normality of $N$, we have $G^{\mathscr{U}_p}\cap N=1$ or $G^{\mathscr{U}_p}\cap N=N$. If the former holds, then $N\leq C_G(G^{\mathscr{U}_p})=1$, a contradiction. If the later holds, then $G^{\mathscr{U}_p}/N\cong M^{\mathscr{U}_p}$. On the other hand, $[N,M^{\mathscr{U}_p}]=1$, so $N\leq C_G(M^{\mathscr{U}_p})\cap C_G(N)=C_G(G^{\mathscr{U}_p})=1$ since $G^{\mathscr{U}_p}=NM^{\mathscr{U}_p}=N\times M^{\mathscr{U}_p}$, a contradiction.  Now if $G=C_G(N)$, then $N\leq C_G(G^{\mathscr{U}_p})=1$, a contradiction. So $C_G(N)<G$ and we may chose a  non-$p$-supersolvable subgroup $L$ of $G$ such that $G=C_G(N)L$ and $C_G(N)\cap L\leq \Phi(L)$. Obviously, $N$ normalize $L^{\mathscr{U}_p}$ and $N\cap L^{\mathscr{U}_p}\unlhd G$. Since  $N$ is a minimal normal subgroup of $G$, we have $N\cap L^{\mathscr{U}_p}=1$ or $N\cap L^{\mathscr{U}_p}=N$. If $N\cap L^{\mathscr{U}_p}=1$, then $L^{\mathscr{U}_p}\leq C_G(N)\cap L$, so $G/C_G(N)\cong L/L\cap C_G(N)\in\mathscr{U}_p$, thus $G^{\mathscr{U}_p}\leq C_G(N)$ and hence $N\leq C_G(G^{\mathscr{U}_p})=1$, a contradiction. Which implies that $N\leq L^{\mathscr{U}_p}$. Now by the fact $G=C_G(N)L$ that $N$ is a minimal normal subgroup of $L$. So $F(L)= C_L(N)=C_G(N)\cap L=\Phi(L)$ and $F(L)<F_p(L)$, as a similar argument of \textbf{Case 1}, we also have a contradiction that $F_p(L)=F(L)$. \qed\\

As  corollaries of Theorem \ref{thmCGNUp1}, we have

\begin{cor}\label{corCGU1}
Let $G$ be a solvable group and $H$ a normal subgroup of $G$, then $N_{\mathscr{U}}(G,H)=1$ if and only if $C_G(H^{\mathscr{U}})=1$.
\end{cor}

\begin{cor}\label{corNUpCGUp1}
Let $G$ be a $p$-solvable group. Then $N_{\mathscr{U}_p}(G)=1$ if and only if $C_G(G^{\mathscr{U}_p})=1$.
\end{cor}

\begin{cor}\label{corZGUpCNUp}
Let $G$ be a $p$-solvable group, if $Z(G^{\mathscr{U}_p})=1$, then $C_G(G^{\mathscr{U}_p})=N_{\mathscr{U}_p}(G)$.
\end{cor}
\proof By Corollary \ref{corNUpCGUp1}, we may assume that $C_G(G^{\mathscr{U}_p})>1$. Moreover, obviously, $C_G(G^{\mathscr{U}_p})\leq N_{\mathscr{U}_p}(G)$. Now  consider $G/C_G(G^{\mathscr{U}_p})$, let $$gC_G(G^{\mathscr{U}_p})\in C_{G/C_G(G^{\mathscr{U}_p})}((G/C_G(G^{\mathscr{U}_p})^{\mathscr{U}_p})=C_{G/C_G(G^{\mathscr{U}_p})}(G^{\mathscr{U}_p}C_G(G^{\mathscr{U}_p})/C_G(G^{\mathscr{U}_p})),$$
 then $[g,G^{\mathscr{U}_p}]\leq C_G(G^{\mathscr{U}_p})\cap G^{\mathscr{U}_p}=Z(G^{\mathscr{U}_p})=1$, thus $g\in C_G(G^{\mathscr{U}_p})$ and $C_{G/C_G(G^{\mathscr{U}_p})}(G/C_G(G^{\mathscr{U}_p})^{\mathscr{U}_p})=\bar{1}$. By Corollary \ref{corNUpCGUp1}, we have $N_{\mathscr{U}_p}(G/C_G(G^{\mathscr{U}_p}))=\bar{1}$, so by Lemma \ref{lemSQ}(2), $N_{\mathscr{U}_p}(G)\leq C_G(G^{\mathscr{U}_p})$, as desired. \qed

\begin{cor}
Let $G$ be a $p$-solvable group, then $Z(G^{\mathscr{U}_p})=1$ if and only if $G^{\mathscr{U}_p}\cap N_{\mathscr{U}_p}(G)=1$.
\end{cor}
\proof Since $Z(G^{\mathscr{U}_p})\leq G^{\mathscr{U}_p}\cap N_{\mathscr{U}_p}(G)$, so the ``only if"  part is obviously. Now assume that $Z(G^{\mathscr{U}_p})=1$, it follows form Corollary \ref{corZGUpCNUp} that $C_G(G^{\mathscr{U}_p})=N_{\mathscr{U}_p}(G)$, so $G^{\mathscr{U}_p}\cap N_{\mathscr{U}_p}(G)=G^{\mathscr{U}_p}\cap C_G(G^{\mathscr{U}_p})=Z(G^{\mathscr{U}_p})=1$, as desired. \qed

\begin{thm}\label{lemZGUpNUp}
Let $G$ be a group, then $Z_{\infty}(G^{\mathscr{U}_p})\leq N^{\infty}_{\mathscr{U}_p}(G)$.
\end{thm}
\proof If $Z(G^{\mathscr{U}_p})=1$, obviously the result holds. So we assume $Z(G^{\mathscr{U}_p})>1$. Now, $Z(G^{\mathscr{U}_p})\leq N_{\mathscr{U}_p}(G)$ by definition. By induction on $|G|$, $Z_{\infty}(G^{\mathscr{U}_p}/Z(G^{\mathscr{U}_p}))=Z_{\infty}((G/Z(G^{\mathscr{U}_p}))^{\mathscr{U}_p})\leq N^{\infty}_{\mathscr{U}_p}(G/Z(G^{\mathscr{U}_p}))$, so $Z_{\infty}(G^{\mathscr{U}_p})\leq N^{\infty}_{\mathscr{U}_p}(G)$ by Lemma \ref{lemDnUp}(3). \qed\\

As a similar argument above, we have

\begin{cor}\label{corZGUNUp}
Let $G$ be a group, then $Z_{\infty}(G^{\mathscr{U}})\leq N^{\infty}_{\mathscr{U}_p}(G)$ for every $p\in\pi(G)$.
\end{cor}

At the end of this section, we investigate the relationship between $N_{\mathscr{U}_p}(G,H)$ and $N_{\mathscr{U}}(G,H)$.

\begin{lem}\label{lemNUSUp}
Let $G$ be a group and $H$ a subgroup of $G$, then $\underset{p\in\pi(H)}{\bigcap}N_{\mathscr{U}_{p}}(G,H)\leq N_{\mathscr{U}}(G,H)$.
\end{lem}
\proof Denote $T=\underset{p\in\pi(G)}{\bigcap}N_{\mathscr{U}_{p}}(G,H)$, then for any subgroup $K$ of $H$, by definition, $T$ normalize $K^{\mathscr{U}_p}$ for every $p\in\pi(K)$. Now by Lemma \ref{lemGUUp}, $T$ normalize $K^{\mathscr{U}}$ for any $K\leq H$, so $T\leq N_{\mathscr{U}}(G,H)$, as desired. \qed

\begin{lem}\label{lemGUisp}
Let $G$ be a group and $H$ a subgroup of $G$. If $H^{\mathscr{U}}$ is a $p$-group for some $p\in\pi(H)$, where $\pi(H)=\{p_1,p_2,\cdots,p_n\}$, then $\underset{p_i\in\pi(H)}{\bigcap}N_{\mathscr{U}_{p_i}}(G,H)= N_{\mathscr{U}}(G,H)$.
\end{lem}
\proof By Lemma \ref{lemNUSUp}, we only need to prove $N_{\mathscr{U}}(G,H)\leq \underset{p_i\in\pi(H)}{\bigcap}N_{\mathscr{U}_{p_i}}(G,H)$. Without loss of generality , we assume $p=p_1$. By Lemma \ref{lemGUUp}, $H^{\mathscr{U}_{p_i}}$ is a $p$-group for every $p_i\in\pi(H)$. Now if $j\neq 1$, the $p_j$-supersolvability of $H/H^{\mathscr{U}_{p_j}}$ implies that $H\in\mathscr{U}_{p_j}$. Thus by Lemma \ref{lemGUUp}, $K^{\mathscr{U}}=K^{\mathscr{U}_p}$ for every subgroup $K$ of $H$ and $G=N_{\mathscr{U}_{p_j}}(G,H)$ for every $p\neq p_j\in\pi(H)$, so $\underset{p_i\in\pi(H)}{\bigcap}N_{\mathscr{U}_{p_i}}(G,H)=N_{\mathscr{U}_{p}}(G,H)=N_{\mathscr{U}}(G,H)$, as desired. \qed

\begin{thm}\label{thmNUpeNU}
Let $G$ be a solvable group and $H$ a normal subgroup of $G$, then $\underset{p\in\pi(H)}{\bigcap}N_{\mathscr{U}_{p}}(G,H)=G$ if and only if $N_{\mathscr{U}}(G,H)=G$.
\end{thm}
\proof Denote $\pi(H)=\{p_1,p_2,\cdots,p_n\}$, the ``$\Rightarrow$" part of this result follows from Lemma \ref{lemNUSUp}. Now we prove the  ``$\Leftarrow$" part, this is equivalent to prove that $N_{\mathscr{U}_{p_i}}(G,H)=G$ for arbitrary $p_i\in\pi(H)$. Assume false, then there is at least one prime, say $p$ in $\pi(H)$ such that $N_{\mathscr{U}_p}(G,H)<G$.  If $p\not\in\pi(H^{\mathscr{U}_p})$, then the $p$-supersolvability of $H/H^{\mathscr{U}_p}$ implies that $H\in\mathscr{U}_p$, then $H^{\mathscr{U}_p}=1$, so $G=N_{\mathscr{U}_p}(G,H)$, a contradiction. By Lemma \ref{lemSQ}(1), $H=N_{\mathscr{U}}(H)$. So  $H^{\mathscr{U}}$ is nilpotent by \cite[Theorem B]{SuNing14}. Note that $p\in\pi(H^{\mathscr{U}})$, so there is a minimal normal subgroup $N$ of $G$ contained in $O_p(H^{\mathscr{U}})$. By Lemma \ref{lemwn}, $G/N=N_{\mathscr{U}}(G/N,H/N)$. By induction on $|G||H|$, we have $G/N=N_{\mathscr{U}_{p_i}}(G/N,H/N)$ for any $p_i\in\pi(H)$. Now by Proposition \ref{proN1}, $(H/N)^{\mathscr{U}_{p_i}}$ is a $p_i$-group for any $p_i\in\pi(H)$. In particular, $H^{\mathscr{U}_p}$ is a $p$-group. If there exists some $p\neq p_i\in\pi(H)$ such that $H^{\mathscr{U}_{p_i}}\neq 1$, then as  a similar proof above,  there is a minimal normal subgroup $T$ of $G$ contained in $O_{p_i}(H^{\mathscr{U}_{p_i}})$ such that $G/T=N_{\mathscr{U}_p}(G/T,H/T)$. For any subgroup $K$ of $H$, we have $(KT/T)^{\mathscr{U}_{p}}\unlhd G/T$, so $K^{\mathscr{U}_{p}}T\unlhd G$. Note that $[K^{\mathscr{U}_{p}},T]\leq [H^{\mathscr{U}_p},T]\leq H^{\mathscr{U}_{p}}\cap T=1$ since $p\neq p_i$, hence $K^{\mathscr{U}_{p}}$ char $K^{\mathscr{U}_{p}}T\unlhd G$, which implies that $K^{\mathscr{U}_{p}}\unlhd G$. This is also contrary to that $G=N_{\mathscr{U}_{p}}(G,H)$. So $H^{\mathscr{U}_{p_i}}=1$ for any $p\neq p_i\in\pi(H)$, thus  by Lemma \ref{lemGUUp}, $H^{\mathscr{U}}=H^{\mathscr{U}_{p_1}}H^{\mathscr{U}_{p_2}}\cdots H^{\mathscr{U}_{p_n}}=H^{\mathscr{U}_p}$ is a $p$-group, now by Lemma \ref{lemGUisp}, we have $\underset{p_i\in\pi(H)}{\bigcap}N_{\mathscr{U}_{p_i}}(G,H)= N_{\mathscr{U}}(G,H)$, a contradiction. \qed

\begin{thm}\label{lemNEU}
Let $G$ be a solvable group and $H$ a normal subgroup of $G$, then $\underset{p\in\pi(H)}{\bigcap}N_{\mathscr{U}_p}(G,H)=1$ if and only if $N_{\mathscr{U}}(G,H)=1$.
\end{thm}
\proof  Assume that $N_{\mathscr{U}}(G,H)=1$, then by Lemma \ref{lemNUSUp}, $\underset{p\in\pi(H)}{\bigcap}N_{\mathscr{U}_{p}}(G,H)=1$. Now assume that $\underset{p\in\pi(H)}{\bigcap}N_{\mathscr{U}_{p}}(G,H)=1$, since $C_G(H^{\mathscr{U}})\leq C_G(H^{\mathscr{U}_p})\leq N_{\mathscr{U}_p}(G,H)$, so $1=C_G(H^{\mathscr{U}})\leq\underset{p\in\pi(H)}{\bigcap}N_{\mathscr{U}_p}(G,H)$. By Corollary \ref{corCGU1}, $N_{\mathscr{U}}(G,H)=1$. \qed

\section{$\mathscr{N}_1$-group}

Let $G$ be a group and $H$ a subgroup of $G$. We now consider the case that $G=N_{\mathscr{U}_p}(G,H)$, i.e.,  $G$ is an $\mathscr{N}_1$-group with respect to $H$ and $\mathscr{U}_p$. Obviously,  $G$ is an $\mathscr{N}_1$-group with respect to $H$  and $\mathscr{U}_p$ if and only if $K^{\mathscr{U}_p}\unlhd G$ for every $K\leq H$. There are many groups which are $\mathscr{N}_1$-groups, here we list some of them.

\begin{pro}
The following groups are $\mathscr{N}_1$-groups with respect to some subgroup $H$ and $\mathscr{U}_p$:

(1) All ($p$-)supersolvable groups.

(2) Groups with  $H$  a ($p$-)supersolvable subgroup.

(3) Groups all of whose non-($p$-)supersolvable subgroups are normal.

(4) Groups with $G^{\mathscr{U}_p}$  a cyclic subgroup.

(5) Groups with  a normal  minimal non-$p$-supersolvable subgroup $H$.

\end{pro}

\begin{pro}\label{proN1G}

Let $G$ be an $\mathscr{N}_1$-group with respect to $H$ and $\mathscr{U}_p$, then

(1) For any $K\leq G$, $K$ is an $\mathscr{N}_1$-group with respect to $H\cap K$ and $\mathscr{U}_p$;

(2) Let $N\unlhd G$ and $N\leq H$, then $G/N$ is an $\mathscr{N}_1$-group with respect to $H/N$ and $\mathscr{U}_p$;

\end{pro}
\proof  The proof of (1) and (2) follows from Lemmas \ref{lemwn}(2)(3) respectively.  \qed\\

In general, if both of a normal subgroup $N$ and corresponding quotient group $G/N$ are $\mathscr{N}_1$-groups with respect to themselves and $\mathscr{U}_p$, then $G$ may be not an $\mathscr{N}_1$-group.

\begin{ex}
Let $G=\langle a,b,c,d\rangle\rtimes\langle e,f\rangle\cong C^4_2\rtimes C_6$, where $e^2=f^3=1$ and $a^e=ac,b^e=bd,c^e=c,d^e=d,a^f=b,b^f=ab,c^f=d,d^f=cd$ (IdGroup=[96,70] in GAP \cite{GAP}). Let $N=\langle a,b,c,d,e\rangle$, then $N=N_{\mathscr{U}_2}(N)$ and $G/N=N_{\mathscr{U}_2}(G/N)$. Now let $H=\langle a,b,f\rangle\cong A_4$, the alternating group of degree $4$, then $H^{\mathscr{U}_2}=\langle a,b\rangle$ is not normal in $G$, so $G$ is not an $\mathscr{N}_1$-group with respect to itself and $\mathscr{U}_2$.
\end{ex}

In the case $\mathscr{F}\in\{\mathscr{A},\mathscr{N}\}$, the solvability of $N_{\mathscr{F}}(G)$ has been investigated by many scholars. For example, see \cite{LiSR10,ShenZC12}. For the case $\mathscr{F}=\mathscr{U}$, we have

\begin{thm}\label{thmNUS}
Let $G$ be an $\mathscr{N}_1$-group with respect to itself and $\mathscr{U}$. Then $G$ is solvable. In particular, $N_{\mathscr{U}}(G)$ and $N^{\infty}_{\mathscr{U}}$ are solvable.
\end{thm}
\proof Obviously, every subgroup and every homomorphic image of $G$ satisfies our hypothesis. If $G$ is not solvable, then $G$ is a non-abelian simple group. Let $H$ be  a proper subgroup of $G$, then $H^{\mathscr{U}}\unlhd G$, so $H^{\mathscr{U}}=1$ and $H$ is supersolvable. This implies that $G$ is a minimal non-supersolvable group. Now by Doerk's result \cite{Doerk66}, $G$ is solvable, a contradiction. \qed\\

Let $H$ be a subgroup of $G$, for the case of $\mathscr{F}=\mathscr{U}_p$, the solvability of $N_{\mathscr{U}_p}(G,H)$  does not always holds. Here, we have

\begin{thm}\label{thmHSpS}
Let $G$ be a group and $H$ a   subgroup of $G$, if $G$ is an $\mathscr{N}_1$-group with respect to $H$  and $\mathscr{U}_p$, where $p$ is the smallest prime in $\pi(H)$, then $H$ is $p$-solvable.

\end{thm}
\proof Assume that the result is false and let the pair $(G,H)$ be a counterexample with $|G||H|$ minimal. By Feit-Thompson's theorem, we  may assume that $p=2$. If $H<G$, then by Lemma \ref{lemwn}(2), $H=N_{\mathscr{U}_2}(H,H)=N_{\mathscr{U}_2}(H)$, so $H$ is $2$-solvable by the choice of $(G,H)$, a contradiction. Now we assume that $H=G$, by hypothesis, $G=N_{\mathscr{U}_2}(G)$.   If $G$ is a  non-abelian simple group, we prove the contradiction  that $G$ is a solvable simple group. Let $K$ be any proper subgroup of $G$, then $K^{\mathscr{U}_2}\unlhd G$, so $K^{\mathscr{U}_2}=1$ and $K$ is $2$-supersolvable. Which implies that $G$ is a non-$2$-supersolvable  group  all of whose proper subgroups are $2$-supersolvable. Now by \cite[1.1]{Tuccillo92}, $G$ is solvable, a contradiction. Let $N$ be a minimal normal subgroup of $G$, then Lemma \ref{lemSQ} implies that $G/N=N_{\mathscr{U}_2}(G/N)$
 and $N=N_{\mathscr{U}_2}(N)$, thus $G/N$ and $N$ are $2$-solvable by our choice of $(G,H)$, hence $G$ is $2$-solvable, the finial contradiction complement our proof.  \qed

\begin{cor}
Let $G$ be a group and $p$ is the smallest prime in $\pi(G)$. If $G$ is an $\mathscr{N}_1$-group with respect to itself and $\mathscr{U}_p$, then  $G$ is $p$-solvable.
\end{cor}

In Theorem \ref{thmHSpS}, the restriction that $p$ is smallest prime in $\pi(H)$ is necessary. Furthermore, in general, we can't obtain the $p$-solvability of $G$ even though $p=2$.

\begin{ex} Let $G=SL(2,5)$(IdGroup=[120,5], see \cite{GAP}), then $\pi(G)=\{2,3,5\}$ and $G$ is not $2$-solvable and not $3$-solvable. But $G=N_{\mathscr{U}_3}(G)$.

Furthermore, $G$ has a normal subgroup, say $H\cong C_2$, obviously, $G=N_{\mathscr{U}_2}(G,H)$. But we have the following obvious conclusions:
\end{ex}

\begin{thm}
Let $G$ be a group and $H$ a subgroup of $G$ . Assume $G$ is an $\mathscr{N}_1$-group with respect to $H$  and $\mathscr{U}_p$, where $p$ is the smallest prime divisor of $|H|$. Then $G$ is $p$-solvable if   one of the statement holds:

(1) $H\unlhd G$ and $G/H$ is $p$-solvable;

(2) $G^{\mathscr{U}_p}\leq H$;




\end{thm}

\begin{thm}
Let $G$ be a group and $P$ a Sylow $p$-subgroup of $G$. Assume that $P\unlhd G$ and $G$ is an $\mathscr{N}_1$-group with itself and $\mathscr{U}_p$, then $G$ is solvable.
\end{thm}

\begin{thm}\label{thmpnhlp}
Let $G$ be a  $p$-solvable  $\mathscr{N}_1$-group with respect to itself  and $\mathscr{U}_p$, then

(1) $G^{\mathscr{U}_p}$ is a $p$-group.

(2) $N_{\mathscr{U}_p}(G/N)>1$ for any proper normal subgroup $N$ of $G$.

(3) $h_p(G)\leq 3$;

\end{thm}
\proof (1) It follows from Proposition \ref{proN1}.

(2) Let $N$ be a proper normal subgroup of $G$, by Lemma \ref{lemSQ}(2), $G/N=N_{\mathscr{U}_p}(G/N)$, so $N_{\mathscr{U}_p}(G/N)>1$ by our choice of $N$.

(3) By (1), $G^{\mathscr{U}_p}$ is a $p$-group, then $G^{\mathscr{U}_p}\leq F_p(G)$.  Consider $\bar{G}=G/F_p(G)$, then $\bar{G}$ is $p$-supersolvable, thus $\bar{G}'=G'F_p(G)/F_p(G)$ is $p$-nilpotent. So $h_p(G)\leq h_p(\bar{G})+1\leq h_p(\bar{G}/\bar{G}')+h_p(\bar{G}')+1=3$. \qed


As  a corollary of Theorem \ref{thmpnhlp}, we have

\begin{cor}\label{GUp32}
Let $G$ be a   $\mathscr{N}_1$-group with respect to itself  and $\mathscr{U}_p$ for every $p\in\pi(G)$, then

(1) $G$ is solvable.

(2) $G=N_{\mathscr{U}}(G)$.

(3) $G^{\mathscr{U}}$ is nilpotent.

(4) $G/N\in\mathscr{N}\mathscr{U}$ for every normal subgroup $N$ of $G$.

(5) $h(G)\leq 3$.

(6) $l_r(G)\leq 2$ for every $r\in\pi(G)$.

\end{cor}
\proof (1) Let $p\in\pi(G)$ be the smallest prime, if $p>2$, obviously $G$ is solvable by the well know Feit-Thompson Theorem. If $p=2$, then by Theorem \ref{thmHSpS}, $G$ is $2$-solvable, so $G$ is solvable, as desired.

(2) By hypothesis, $G=N_{\mathscr{U}_p}(G)$ for every $p\in\pi(G)$.  Now by Theorem \ref{thmNUpeNU}, $G=N_{\mathscr{U}}(G)$.

(3) Let $p\in\pi(G)$, if $G\in\mathscr{U}_p$, then $G^{\mathscr{U}_p}=1$, if $G\not\in\mathscr{U}_p$, by Theorem \ref{thmpnhlp}(1), $G^{\mathscr{U}_p}$ is a $p$-group. Now assume $\pi(G)=\{p_1,p_2,\cdots,p_n\}$, then by  Lemma \ref{lemGUUp}, $G^{\mathscr{U}}=G^{\mathscr{U}_{p_1}}G^{\mathscr{U}_{p_2}}\cdots G^{\mathscr{U}_{p_n}}$, so $G^{\mathscr{U}}$ is nilpotent.

(4) By (2), $G^{\mathscr{U}}$ is nilpotent. Let $N$ be a proper normal subgroup of $G$, then $(G/N)^{\mathscr{U}}=G^{\mathscr{U}}N/N\cong G^{\mathscr{U}}/G^{\mathscr{U}}\cap N$ is nilpotent, as desired.

(5) It follows from Theorem \ref{thmpnhlp}(3).

(6) If $O_{r'}(G)>1$, then $N_{\mathscr{U}_p}(G/O_{r'}(G))=G/O_{r'}(G)$ for every $p\in\pi(G)$ by Lemma  \ref{lemSQ}(2), so $G/O_{r'}(G)$ satisfies our hypothesis. By induction, $l_r(G/O_{r'}(G))\leq 2$, which implies that $l_r(G)\leq 2$. As a similarly argument above and by Lemma \ref{lempL}, we have $l_r(G)=l_r(G/\Phi(G))\leq 2$, hence we may assume that $O_{r'}(G)=\Phi(G)=1$, so $F(G)=F_r(G)=O_r(G)$. Now consider $G/F(G)$, by (5), $h(G/F(G))\leq 2$. Now by Lemma \ref{lempl21}, $l_r(G/F(G))=l_r(G/O_r(G))\leq 1$, so $l_r(G)\leq l_r(G/O_r(G))+1\leq 2$, as desired.  \qed\\
%

Recall that the well-know results of P.Hall and D.J.S. Robinson  as follows:

\begin{thm} Let $G$ be a group and $N$ a nilpotent normal subgroup of $G$. Then

(1) (P. Hall)  If $G/N'$ is nilpotent, then $G$ is nilpotent.

(2) (D.J.S. Robinson). If $G/N'$ is supersolvable, then $G$ is supersolvable.
\end{thm}

Here we have

\begin{thm}
Let $G$ be a $p$-solvable $\mathscr{N}_1$-group with respect to itself and $\mathscr{U}_p$. Then $G\in\mathscr{N}$  if $G^{\mathscr{N}}\leq G^{\mathscr{U}^2_p}$.
\end{thm}

\begin{thm}
Let $G$ be an $\mathscr{N}_1$-group with respect to itself and $\mathscr{U}$. Then $G\in\mathscr{U}$  if $G^{\mathscr{U}}\leq G^{\mathscr{U}^2}$.
\end{thm}

\section{Applications}

\begin{thm}\label{thme2}
Let $G$ be a $p$-solvable group, then the following statements are equivalent:

(1) $G\in\mathscr{N}_p\mathscr{U}_p$.

(2) $G/N_{\mathscr{U}_p}(G)\in\mathscr{N}_p\mathscr{U}_p$.
\end{thm}
\proof (1)$\Rightarrow$ (2). It is obviously.

(2) $\Rightarrow$ (1). Assume the result is false and let $G$ be a counterexample of minimal order. If $N_{\mathscr{U}_p}(G)=1$, there is nothing to prove. Now we assume that $N_{\mathscr{U}_p}(G)>1$. Let $N$ be a minimal normal subgroup of $G$ contained in $N_{\mathscr{U}_p}(G)$. Then either $N$ is an elementary abelian $p$-group or a $p'$-group. Note that $G/N/N_{\mathscr{U}_p}(G/N)\leq G/N/N_{\mathscr{U}_p}(G)/N\cong G/N_{\mathscr{U}_p}(G)$ by Lemma \ref{lemSQ}(2), then $G/N\in\mathscr{N}_p\mathscr{U}_p$ by the choice of $G$. So we may assume that $N$ is an elementary abelian $p$-group, otherwise, $G^{\mathscr{U}_p}N/N\cong G^{\mathscr{U}_p}/(G^{\mathscr{U}_p}\cap N) \in\mathscr{N}_p$, so $G\in\mathscr{N}_p\mathscr{U}_p$ by Lemma \ref{lempnPhi}(2), a contradiction. Furthermore,  if $N\leq\Phi(G)$, then $G/\Phi(G)\in\mathscr{N}_p\mathscr{U}_p$, so $G^{\mathscr{U}_p}\Phi(G)/\Phi(G)\in\mathscr{N}_p$, hence $G\in\mathscr{N}_p\mathscr{U}_p$ by Lemma \ref{lempnPhi}(1), a contradiction again. Now there exist some maximal subgroup $M$ of $G$ such that $G=MN$ and $M\cap N=1$. Also we have $M^{\mathscr{U}_p}=(G/N)^{\mathscr{U}_p}\in\mathscr{N}_p$. Note that $N\leq N_{\mathscr{U}_p}(G)$, then $N$ normalize $M^{\mathscr{U}_p}$, thus $[N,M^{\mathscr{U}_p}]=1$ and $NM^{\mathscr{U}_p}=N\times M^{\mathscr{U}_p}\unlhd G$. Since $G/NM^{\mathscr{U}_p}=MN/NM^{\mathscr{U}_p}=MNM^{\mathscr{U}_p}/NM^{\mathscr{U}_p}\cong M/M\cap NM^{\mathscr{U}_p}=M/M^{\mathscr{U}_p}\in\mathscr{U}_p$, then $G^{\mathscr{U}_p}\leq NM^{\mathscr{U}_p}\in\mathscr{N}_p$, a contradiction. \qed

\begin{cor}\label{GEDnUp}
Let $G$ be a $p$-solvable group, if $G=N^{\infty}_{\mathscr{U}_p}(G)$, then

(1) $G\in\mathscr{N}_p\mathscr{U}_p$.

(2) $N^{\infty}_{\mathscr{U}_p}(G/N)>1$ for any proper normal subgroup $N$ of $G$.

(3) $h_p(G)\leq 3$;

\end{cor}
\proof (1) Set $G=N^n_{\mathscr{U}_p}(G)$ for some positive integer $n$. If $n=1$, the result follows from Theorem \ref{thmpnhlp}(1). If $n\geq 2$, then $N_{\mathscr{U}_p}(G)<G$. By Lemma \ref{lemDnUp}, $G/N_{\mathscr{U}_p}(G)=N^{\infty}_{\mathscr{U}_p}(G/N_{\mathscr{U}_p}(G))$, so $G/N_{\mathscr{U}_p}(G)$ satisfies our hypothesis, thus $G/N_{\mathscr{U}_p}(G)\in\mathscr{N}_p\mathscr{U}_p$ by induction on $|G|$. Now by Theorem \ref{thme2}, $G\in\mathscr{N}_p\mathscr{U}_p$, as desired.

(2) If there exists some proper normal subgroup $N$ of $G$ such that $N^{\infty}_{\mathscr{U}_p}(G/N)=1$, then $G=N^{\infty}_{\mathscr{U}_p}(G)\leq N$ by Lemma \ref{lemDnUp}(3), that is impossible.

The proof of (3) is similar to Theorem \ref{thmpnhlp}(3). \qed\\

We have proved in Theorem \ref{thmNUS} that $G$ is solvable if $G=N^{\infty}_{\mathscr{U}}(G)$. Now as similar argument above, we have

\begin{cor}\label{croShlp}
Let $G$ be a group, if $G=N^{\infty}_{\mathscr{U}}(G)$, then

(1) $G\in\mathscr{N}\mathscr{U}$.


(2) $h(G)\leq 3$.

(3) $l_p(G)\leq 2$ for every $p\in\pi(G)$.
\end{cor}

\begin{thm}\label{thmlhp}
Let $G$ be a group and $p\in\pi(G)$.

(1) Assume $G$ is $p$-solvable, then $h_p(G)\leq k$ if and only if  $h_p(G/N^{\infty}_{\mathscr{U}_p}(G))\leq k$, where $k\geq 3$.

(2) Assume $G$ is solvable, then $l_p(G)\leq k$ if and only if  $l_p(G/N^{\infty}_{\mathscr{U}}(G))\leq k$, where  $k\geq 2$.
\end{thm}
\proof   We only prove (1) since the proof of (2) is similar and based on  Corollary \ref{croShlp}(4). The ``$\Rightarrow$" part of this theorem is obviously, so we prove the``$\Leftarrow$" part. Setting $N^{\infty}_{\mathscr{U}_p}(G)=N^n_{\mathscr{U}_p}(G)$ for some positive integer $n$.   If $G=N^n_{\mathscr{U}_p}(G)$, then the result follows form Corollary \ref{GEDnUp}(3). So we assume that $N^n_{\mathscr{U}_p}(G)<G$.

Assume that the result is false and let the pair $(G,n)$ be a counterexample with $n|G|$ minimal. For the case $n=1$, we have $h_p(G/N_{\mathscr{U}_p}(G))\leq k$ by hypothesis. Let $N$ be a minimal normal subgroup of $G$, then by Lemma \ref{lemSQ}(2), $G/N/N_{\mathscr{U}_p}(G/N)\leq G/N/N_{\mathscr{U}_p}(G)N/N\cong G/N_{\mathscr{U}_p}(G)N\cong G/N_{\mathscr{U}_p}(G)/N_{\mathscr{U}_p}(G)N/N_{\mathscr{U}_p}(G)$, so $h_p(G/N/N_{\mathscr{U}_p}(G/N))\leq k$ by Lemma \ref{lemhp}(2), thus $h_p(G/N)\leq k$ by the choice of $(G,1)$. Suppose that $N$ is not unique, that is $G$ has at least two minimal normal subgroups, say $N_1$, $N_2$ and $N_1\neq N_2$, then as above, $h_p(G/N_1),h_p(G/N_2)\leq k$, so $h_p(G)\leq k$ by Lemma \ref{lemhp}(4), a contradiction. With similar argument and by using of Lemma \ref{lemhp}(5), we have $\Phi(G)=1$. Furthermore, if $O_{p'}(G)>1$, then it follows  from Lemma \ref{lemhp}(6) that $h_p(G/O_{p'}(G))=h_p(G)\leq k$, a contradiction. So  $F_p(G)=F(G)=C_G(N)=N\leq N_{\mathscr{U}_p}(G)$. Therefore, $G=MN$, $M\cap N=1$ for some maximal subgroup $M$ of $G$. If $M^{\mathscr{U}_p}=1$, that is $G/N$ is $p$-supersolvable, so $G/N=N_{\mathscr{U}_p}(G/N)=N^{\infty}_{\mathscr{U}_p}(G/N)=N^{\infty}_{\mathscr{U}_p}(G)/N$, thus $G=N^{\infty}_{\mathscr{U}_p}(G)$, a contradiction.  If $M^{\mathscr{U}_p}>1$, then $N$ normalize $M^{\mathscr{U}_p}$, so $[N,M^{\mathscr{U}_p}]=1$ and hence $M^{\mathscr{U}_p}\leq C_G(N)=N$, $M^{\mathscr{U}_p}\leq N\cap M=1$, a contradiction again. Now assume that $n>1$, since $G/N^n_{\mathscr{U}_p}(G)\cong G/N^{n-1}_{\mathscr{U}_p}(G)/N^n_{\mathscr{U}_p}(G)/N^{n-1}_{\mathscr{U}_p}(G)=G/N^{n-1}_{\mathscr{U}_p}(G)/N_{\mathscr{U}_p}(G/N^{n-1}_{\mathscr{U}_p}(G))$, so by induction and combining the proof of the case $n=1$, we have that  $h_p(G/N^{n-1}_{\mathscr{U}_p}(G))\leq k$ , hence  $h_p(G)\leq k$ by our choice of $(G,n)$, a contradiction. \qed

\begin{thm}\label{thmlh23}
Let $G$ be a $p$-solvable group. If all elements of  order $p$ and order $4$(if $p=2$) of $G$ are in $N^{\infty}_{\mathscr{U}_p}(G)$, then $h_p(G)\leq 3$.
%
%
\end{thm}
\proof  Assume that the result is false and let $G$ be a counterexample with minimal order, we prove this result by several steps:

\textbf{Step 1}. Let $H<G$, then $h_p(H)\leq 3$.

Let $H$ be a proper subgroup of $G$, then all elements of $H$ of order $p$ and order $4$(if $p=2$) are contained in $H\cap N^{\infty}_{\mathscr{U}_p}(G)\leq N^{\infty}_{\mathscr{U}_p}(H)$ by Lemma \ref{lemDnUp}(1), so $H$  satisfies our hypothesis, $h_p(H)\leq 3$ by the choice of $G$.

\textbf{Step 2}. $O_{p'}(G)=1$, $C_G(O_p(G))\leq O_p(G)=F(G)$.

Let $T=O_{p'}(G)$, if $T\neq 1$, for any $xT\in N^{\infty}_{\mathscr{U}_p}(G)T/T$ with $|x|=p$ or $|x|=4$(if $p=2$),   $xT\leq N^{\infty}_{\mathscr{U}_p}(G/T)$ by Lemma \ref{lemDnUp}(3), so $G/T$ satisfies our hypothesis, thus $h_p(G/T)\leq 3$ by the choice of $G$. It is easy to see that $h_p(G)\leq 3$, a contradiction. Further more, we have $C_G(O_p(G))\leq O_p(G)=F(G)$ since $G$ is $p$-solvable and $O_{p'}(G)=1$.

\textbf{Step 3}. $\Phi(G)<O_p(G)$.

By (2), $\Phi(G)\leq O_p(G)$, if $O_p(G)=\Phi(G)$, then $O_{pp'}(G)/\Phi(G)=O_{p'}(G/\Phi(G))>1$ since $G$ is $p$-solvable and $O_p(G/\Phi(G))=1$. Let $O_{pp'}(G)=\Phi(G)T$, where $T>1$ is a $p'$-group, then $G=N_G(T)O_{pp'}(G)=N_G(T)$ by Frattini argument, it follows that $T\unlhd G$ and then $T\leq O_{p'}(G)=1$ by \textbf{Step} 2, a contradiction.

\textbf{Step 4}. Finial contradiction.

Denote by $\mathscr{F}_p$ the class of all $p$-solvable groups whose $p$-Fitting length are less than $3$, then by Lemma \ref{lemhp}, $\mathscr{F}_p$ is a saturated formation. Now by \textbf{Step} 1, $G$ is an $\mathscr{F}_p$-critical group, i.e, $G$ is not belongs to $\mathscr{F}_p$, but all proper subgroups  of $G$  belong to $\mathscr{F}_p$. Now by \textbf{Steps} 1,3, there exists some maximal subgroup $M$ of $G$ such that $G=O_p(G)M=F(G)M$ and $M\in\mathscr{F}_p$, so by Lemma \ref{lemAdolf}(2), $G^{\mathscr{F}_p}$ is a $p$-group and $G^{\mathscr{F}_p}$ has exponent $p$ if $p>2$ and exponent at most $4$ if $p=2$, hence $G^{\mathscr{F}_p}\leq O_p(G)$ and $G^{\mathscr{F}_p}\leq N^{\infty}_{\mathscr{U}_p}(G)$, it follows that $h_p(G/N^{\infty}_{\mathscr{U}_p}(G))\leq 3$. Now by Theorem \ref{thmlhp}(2), $h_p(G)\leq 3$, a contradiction. \qed

\begin{cor}\label{corlh23}
Let $G$ be a $p$-solvable group, if every cyclic subgroup of order prime $p$ and order $4$(if $p=2$)  of $G$  is contained in $N_{\mathscr{U}_p}(G)$, then $h_p(G)\leq 3$.
\end{cor}

Combining Theorem \ref{thmlhp}, Lemma \ref{lempL} and Lemma \ref{lemAdolf}, we have the following Theorem \ref{thmlp2}. It's proof is similar to Theorem \ref{thmlh23}.

\begin{thm}\label{thmlp2}
Let $G$ be a solvable group, if every cyclic subgroup of order prime $p$ and order $4$(if $p=2$)  of $G$  is contained in $N^{\infty}_{\mathscr{U}}(G)$, then $l_p(G)\leq 2$.
\end{thm}

As a local vision of Theorem \ref{thmlp2}, we have the following question.

\begin{ques}\label{queslp2}
Let $G$ be a $p$-solvable group, if every cyclic subgroup of order prime $p$ and order $4$(if $p=2$)  of $G$  is contained in $N^{\infty}_{\mathscr{U}_p}(G)$, does $l_p(G)\leq 2$?
\end{ques}

\begin{remark}

(1) In Theorem \ref{thmlh23}, Corollary \ref{corlh23} and Question \ref{queslp2}, the Theorem \ref{thmHSpS} implies the hypothesises that $G$ is $p$-solvable are necessary. For example, consider $G=A_5$, the alternating group of degree $5$, then $G=N_{\mathscr{U}_3}(G)=N_{\mathscr{U}_5}(G)$, so $G$ satisfies our condition, but $G$ is simple.

(2) In Theorem \ref{thmlp2}, the solvability of $G$ is necessary. For example, choice the non-solvable group $G=C_7\times A_5$(IdGroup=[420,13])  and let  $p=7$. Then $N_{\mathscr{U}}(G)=N^{\infty}_{\mathscr{U}}(G)=C_7$ and every cyclic subgroup of order $7$ is in $C_7$.

(3) In Question \ref{queslp2}, the integer 2 is the minimum upper bound. For example, let $G=C_5\times S_4$(IdGroup=[120,37] in GAP). Obviously, $G$ is $2$-solvable and the proper non-2-supersolvable groups of $G$ isomorphism to one of $A_4,S_4,C_5\times A_4$. As $A_4,S_4,C_5\times A_4$ are normal in $G$, so $G=N_{\mathscr{U}_2}(G)$. Note that  $1<C_5<C_{10}\times C_2<C_5\times A_4<C_5\times S_4=G$ is the upper $2$-series of $G$, so $l_2(G)=2$.

\end{remark}


\end{CJK}
\end{document}